\documentclass[12pt,a4paper]{article}

\usepackage{amsmath,amssymb,amsthm}
\usepackage[all]{xy}

\DeclareMathOperator{\Str}{Str}
\DeclareMathOperator{\ann}{Ann}

\newcommand{\zfour}{\mathbb Z_{4}}
\newcommand{\ztwo}{\mathbb Z_{2}}
\newcommand{\zn}{\mathbb Z_{n}}

\newtheorem{unnumber}{}

\newtheorem{prepf}[unnumber]{Proof}
\newenvironment{pf}{\prepf\rm}{\endprepf}

\newcommand{\lcur}{\ar@/^/}
\newcommand{\zfgraph}{\xymatrix@-1.2pc @ur{ 3\lcur@{-}[r] & 2\lcur@{-}[d]
\\ 0\lcur@{-}[u] & 1\lcur@{-}[l]  } }
\newcommand{\ztgraph}{\xymatrix@-1.2pc @ur { 10\lcur@{-}[r] & 11\lcur@{-}[d]
\\ 00\lcur@{-}[u] & 01\lcur@{-}[l] } }

\newtheorem{lem}{Lemma}[section]
\newtheorem{prop}[lem]{Proposition}
\newtheorem{thm}[lem]{Theorem}
\newtheorem{cor}[lem]{Corollary}

\begin{document}

\title{$\zfour$-codes and their Gray map images as orthogonal arrays}
\author{ Peter J Cameron,\footnote{Current address: School of Mathematics and
Statistics, North Haugh, St Andrews, Fife KY16 9SS, UK}
\, Josephine Kusuma
 \\ \small{School of Mathematical Sciences, Queen Mary University of London,} \\
\small{Mile End Road, London E1 4NS, UK}\\ and \\
Patrick Sol\'e \\
\small{Telecom ParisTech, 46 rue Barrault, 75634 Paris Cedex 13, France}\\
\small{\&  Math Dept, King Abdulaziz University, Jeddah, Saudi Arabia}}

\date{}

\maketitle

\begin{abstract}
A classic result of Delsarte connects the strength (as orthogonal array)
of a linear code with the minimum weight of its dual: the former is one
less than the latter.

Since the paper of Hammons \emph{et al.}, there is a lot of interest in
codes over rings, especially in codes over $\zfour$ and their (usually
non-linear) binary Gray map images.

We show that Delsarte's observation extends to codes over arbitrary
finite commutative rings with identity. Also, we show that the strength of
the Gray map image of a $\zfour$ code is one less than the minimum Lee weight
of its Gray map image.
\end{abstract}

%Section 1
\section{Introduction}
Delsarte~\cite{del}, see also \cite[pp. 54--55]{hed},
 showed that the strength (as orthogonal array) of a
linear code over a field is one less than the minimum weight of the dual
code. Two natural questions are: what happens for codes over rings? and
what happens for non-linear codes?

The first question can be answered completely. We show that a linear
code over any finite commutative ring with identity has strength one
less than the minimum weight of its dual (Theorem~\ref{th1}).

The two questions are related. In a groundbreaking paper, Hammons
\emph{et~al.}~\cite{ham} showed that certain non-linear codes which
have more codewords (and thus better error-correcting capabilities)
than any linear code are images of linear codes over $\zfour$ under
the non-linear \emph{Gray map}. (These include the Nordstrom--Robinson,
Preparata and Kerdock codes.) So we are led to investigate the
relationship between the strength of a linear $\zfour$-code and the
strength of its Gray map image. We show
that the strength of the Gray map image of a linear $\zfour$-code
is one less than the minimum Lee weight of the dual $\zfour$-code.
As a consequence, if a linear $\zfour$-code $C$ has strength $t$, then the
strength of its Gray map image lies between $t$ and $2t-1$.

%Section 2
\section{Definitions and preliminaries}

We will begin with some background material on rings and modules, codes over
$\zfour$, the Gray map and orthogonal arrays.

\subsection{Rings and modules}

All our rings will be finite commutative rings with identity, and our
modules are unital. A (linear) code over a ring $R$ is a submodule of
the free module $R^n$; the dual of $C$ is, as usual,
\[C^\bot=\{(x_1,\ldots,x_n)\in R^n:(\forall c_1,\ldots,c_n\in C)\,\,
x_1c_1+\cdots+x_nc_n=0\}.\]
The definitions of zero-divisor, unit, and ideal in a ring, and annihilator
$\ann(M)$ of a module, are standard: see~\cite{hun,mat}.

A ring is \emph{local} if it has a unique maximal ideal, and \emph{semilocal}
if it has only finitely many maximal ideals. If $I$ is an ideal of a ring
$R$, the \emph{$I$-adic completion} of $R$ is the inverse limit of the
rings $R/I^n$ over natural numbers $n$, relative to the natural projections
from $R/I^n$ to $R/I^m$ for $m\le n$. We use the following result
(\cite[p.~62]{mat}):

\begin{thm}
Let $R$ be a semilocal ring with maximal ideals $M_1,\ldots,M_r$, and set
$I=M_1\cdots M_r$ (the radical of $R$). Then the $I$-adic completion of $R$
decomposes as
\[\hat{R}=\hat{R_1}\times\cdots\times\hat{R_n},\]
where $R_i$ is the localisation of $R$ at $M_i$, and $\hat{R_i}$ is its
$M_i$-adic completion.
\label{t:semilocal}
\end{thm}

Now let $R$ be any finite commutative ring with identity. Clearly $R$ is
semilocal, and satisfies the descending chain condition on ideals; so its
radical $I$ is nilpotent, say $I^k=\{0\}$. (See \cite[p.~430]{hun}.) Thus the
$I$-adic completion of $R$ is equal to $R/I^k = R$. Thus, every such ring $R$
is a direct sum of local rings.

\subsection{Codes}

A code $C$ of length $n$ over a ring $R$ (commutative with identity)
is a submodule of the free module $R^n$. Just as in the classical
theory, codewords are $n$-tuples from $R$, and addition and scalar
multiplication are coordinatewise. We define the \textit{Hamming distance}
$d_H(x,y)$, where $x$ and $y$ are the codewords in $C$, as the
number of coordinate places in which $x$ and $y$ differ. The
\textit{Hamming weight} $w_H$ of a codeword $x$ is defined to be
the number of non-zero entries of $x$, $w_H(x) = d_H(x,0)$, where
$0$ is the all zero word in $C$. The minimum Hamming distance of
$C$ is equal to the minimum Hamming weight
among all non-zero codewords in $C$, which we denote by $w_H(C)$.

\subsection{Codes over $\zfour$}

The four elements of $\zfour$ are arranged around a circle. We
define the \textit{Lee distance} $d_L$ as the number of steps
apart they are between each other:
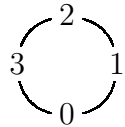
\[\zfgraph\]
This is extended to codewords over $\zfour$ by the rule
$d_L(a,b)=\sum_i d_L(a_i,b_i)$. Similarly, the \textit{Lee weight}
of a codeword $x$ in $C$ is $w_L(x) = d_L(x,0)$.
Moreover, the minimum Lee distance of $C$ is equal
to the minimum Lee weight among all non-zero codewords in $C$; this
is denoted $w_L(C)$.

The \textit{Gray map}, takes $\zfour$ to $\ztwo^2$ by
\[0\mapsto00,\quad1\mapsto01,\quad 2\mapsto11,\quad3\mapsto10,\]
and maps $\zfour^n$ to $\ztwo^{2n}$ coordinatewise.
We denote the Gray map image of a linear code $C$ as $C'$. Note that $C'$ is
usually non-linear.

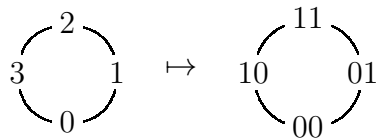
\begin{center}
\begin{tabular} {lcr}
\zfgraph & $\mapsto $ & \ztgraph
\end{tabular}
\end{center}

The Lee weight of $v \in C$
is the Hamming weight of $v' \in C'$, also $d_L(v,w) =
d_H(v',w')$, where $v$, $w \in C$ and $v'$, $w' \in C'$.
In other words, the Gray map, though non-linear, is an isometry from
$\zfour^n$ to $\ztwo^{2n}$.

\subsection{Orthogonal arrays}

An $N \times k$ array $A$ with entries from an alphabet $F$ is
said to be an \textit{orthogonal array} with $f$ levels, where
$|F| = f$, strength $t$ and index $\lambda$ (for some $t$ in the range
$0 \leq t \leq k$) if every $N \times t$ subarray of $A$ contains
each $t$-tuple based on $F$ exactly $\lambda$ times as a row. We say
that a code $D$ is an \textit{orthogonal array} of strength $t$
and index $\lambda$ if the codewords of $D$ make up the rows of an
orthogonal array $A$ as described above. $D$ may or may not be
linear.

Clearly, if $D$ is an orthogonal array of strength $t$ then it is also an
orthogonal array of strength $s$ for any $s \leq t$; its index is given
by $\lambda_s=|F|^{t-s}\lambda$.

%Section 3
\section{Codes over rings as orthogonal arrays}

In this section, we prove that Delsarte's theorem extends to codes over rings.

\begin{thm}
Let $C$ be a linear code over a finite commutative ring $R$ with identity.
Then $\Str(C)=w_H(C^\bot)-1$.
\label{th1}
\end{thm}

\begin{pf}
The proof proceeds in several steps. First, we define two properties which
a ring may have:
\begin{itemize}\itemsep0pt
\item $R$ has property $(+)$ if, whenever $I$ is an ideal of $R$ with
$I\ne R$, we have $\ann(I)\ne0$.
\item $R$ has property $(*)_t$ if given any $N
\times t$ matrix $A$ over ring $R$ with linearly independent
columns, the row space of $A$, $\rho(A)$, is equal to $R^t$.
\end{itemize}

Now the conclusion of the main theorem holds if
$R$ has property $(*)_t$. For suppose that the minimum weight of
$C^\perp$ is $d$. Let $M_C$ be the $n\times N$ matrix whose rows are all
the codewords in $C$, where $N=|C|$. By assumption, any $d-1$ columns of
$M_C$ are linearly independent, since a linear dependence relation would
give a word of weight at most $d-1$ in $C^\perp$. Property $(*)_{d-1}$
shows that the row space of the submatrix formed by these columns is $R^{d-1}$,
so that any $(d-1)$-tuple occurs in these positions; and the number of times
it occurs is equal to the number of codewords with zeros in these positions,
since these words form a coset of the submodule consisting of such codewords.
Thus $M_C$ is an orthogonal array of strength $d-1$.

Next we show, by induction on $t$, that property $(+)$ implies property
$(*)_t$ for all $t$.

To start the induction, note that property $(*)_1$
asserts that if $a$ is an element of $R$ with $\ann(a)=\{0\}$, then
the ideal generated by $a$ is the whole of $R$; this is the contrapositive
of Property $(+)$ for this ideal.

Now, suppose that $(+)$ implies $(*)_{t-1}$ and let $R$ be a ring
satisfying $(+)$ and $A$ an $N \times t$ matrix over $R$ with
linearly independent columns. By the induction hypothesis, if
$\check{A}$ is the matrix obtained from $A$ by deleting the last
column, then $\rho(\check{A}) = R^{t-1}$. This means that, for $1
\leq i \leq t-1$, $\rho(A)$ contains $(0, \ldots, 1, 0, \ldots, 0,
a_i)$ where $1$ is in position $i$ and $a_i \in R$. Let $I = \{r
\in R: (0, \ldots, 0, r) \in \rho(A)\}$. Then, it is easy to check
that $I$ is an ideal of $R$. If it is equal to $R$, then clearly
$\rho(A) = R^t$. So suppose not. Then there is a non-zero element
$b$ annihilating $I$. We claim that
\begin{equation*}
\rho(A) \subseteq \{(x_1, \ldots, x_t) \in R^t: b(x_t -
\sum^{t-1}_{i=1} a_ix_i) = 0\}.
\end{equation*}
For take $(x_1, \ldots, x_t) \in \rho(A)$. I know that $(0,
\ldots, 0, x_i, 0, \ldots, 0, a_ix_i) \in \rho(A)$, for $i=1,
\ldots, t-1$. I can then subtract the sum of $(0, \ldots, 0, x_i,
0, \ldots, 0, a_ix_i)$, for all $i=1, \ldots, t-1$, from $(x_1,
\ldots, x_t)$, and I get $(0, \ldots, 0, x_t - \sum^{t-1}_{i=1}
a_ix_i) \in \rho(A)$; so $(x_t - \sum^{t-1}_{i=1} a_ix_i) \in I$,
whence $b(x_t - \sum^{t-1}_{i=1} a_ix_i) = 0$. Thus, the vectors
in $\rho(A)$ satisfy a non-trivial relation. In particular, this
relation holds for the entries in each row of $A$; that is, it
holds for the columns of $A$. So the columns are not linearly
independent, contrary to assumption.

Our task now is to show that every finite ring $R$ (commutative with
identity) has property $(+)$. We observed after Theorem~\ref{t:semilocal}
that $R$ is a direct product of finite local rings.

Suppose first that $R$ is a local ring, with unique maximal ideal $J$. Again
by finiteness, every element of $J$ is nilpotent, and so there exists an
integer $k$ such that $J^k=\{0\}$ (and we can assume that $k$ is minimal
with this property). Take $a\in J^{k-1}$ with $a\ne0$. Then $a\in\ann(J)$.
Since every proper ideal $I$ is contained in $J$, we have $a\in\ann(I)$,
as required.

Now suppose that $R=R_1\times\cdots\times R_r$, where $R_i$ has property $(+)$
for $i=1,\ldots,r$. Let $I$ be a proper ideal of $R$. Then
$I=I_1\times\cdots\times I_r$, where $I_i$ is an ideal of $R_i$, and is a
proper ideal for at least one value of $i$, say $i=j$. By $(+)$, $R_j$
contains a non-zero element $a_j$ which annihilates $I_j$. Then
$(0,\ldots,0,a_j,0,\ldots,0)$ annihilates $I$.

The theorem is proved.\qed
\end{pf}

We conclude this section with a similar result.

\begin{prop} \label{dim}
Let $R$ be a finite commutative ring with identity, in which
$|\ann(I)|=|R|/|I|$ holds for any ideal $I$. Then, for any code $C$ of
length $n$ over $R$, we have $|C^\perp|=|R|^n/|C|$.
\label{p:dual}
\end{prop}

We omit the proof, which is a straightforward induction. Note that many
finite rings (for example, finite fields, rings $\zn$) have this property.
An example of a ring not satisfying the property is the ring generated
by $\zfour$ and an element $a$ satisfying $2a=a^2=0$. Then $|R|=8$;
the set $\{0,2,a,a+2\}$ is an ideal which is equal to its annihilator.

%Section 4
\section{Orthogonal arrays and the Gray map}

If $D$ is an orthogonal array, then we denote the strength of
$D$ as $\Str(D)$.

\subsection{The main theorem}

The main result is that the strength of the Gray map image of a linear
$\zfour$-code $C$ is one less than the minimum weight of $C^\perp$.

\begin{thm} \label{OAmain}
Let $C$ be a linear $\zfour$-code and $C'$ its Gray map image. Then
$\Str(C')=w_L(C^\bot)-1$, where $w_L$ denotes the minimum Lee weight.
\end{thm}

\begin{pf}
The code $C'$ is the Gray map image of a $\zfour$-linear code, and hence is
distance invariant~\cite{del}. Thus, its weight enumerator is equal
to its distance enumerator, and the \emph{dual distance} is defined; this
is the index of the first non-vanishing term in the MacWilliams transform
of the weight enumerator of $C'$. Delsarte's results show that the strength
of $C'$ is one less than its dual distance.

On the other hand, the Hamming weight enumerator of $C'$ is equal to the
Lee weight enumerator of $C$. Since the MacWilliams relations hold for the
Lee weight enumerators of a linear $\zfour$ code and its dual~\cite{ham},
the dual distance of $C'$ is equal to the minimum weight of $C^\perp$,
and the theorem is proved.\qed
\end{pf}

\subsection{A consequence}

As a corollary, we have an inequality connecting the strength of a linear
$\zfour$-code and that of its Gray map image.

\begin{cor} \label{OA1}
Let $C$ be a $\zfour$-code and $C'$ its Gray map image. If
$\Str(C)=t$, then $t\leq \Str(C') \leq 2t+1$.
\end{cor}

\begin{pf}
If $\Str(C)=t$, then $w_H(C^\bot)=t+1$. Hence $t+1\le w_L(C^\bot)\le 2t+2$,
since any non-zero coordinate of a word in $C^\bot$ contributes either $1$
or $2$ to its Lee weight. The result follows since $\Str(C')=w_L(C^\bot)-1$
by the main theorem. \qed
\end{pf}

This corollary leads us to the following question: Given $t,t'$ such
that $t \leq t' \leq 2t+1$, can we find a code $C$ such that
$\Str(C)=t$ and $\Str(C')=s$?

%Section 5
\section{Further problems}

Connections such as we have explored here probably do not exist between codes
over $\zn$ or more general rings. Although the Lee distance is
defined over $\zn$, for any $n$, there is no analogue of the Gray map.

Another possible topic would be to find connections between these
codes as orthogonal arrays and the $t$-designs formed by the
supports of their codewords.

We also have the problem of deciding which pairs $(t,t')$ are representable
as the strengths of a $\zfour$-code and its Gray map image, as discussed
after Corollary~\ref{OA1}.


\begin{thebibliography}{9}

\bibitem{cam} Cameron, P.J., Van Lint, J.H. \emph{Design, Graphs, Codes and their Links.} London Mathematical Society Student Texts 22, Cambridge University Press, 1991.

\bibitem{del} Delsarte, P. \emph{Four Fundamental Parameters of a Code and their Combinatorial Significance .} Information and Control 23(5), 1973, pp 407-438.

%\bibitem{fuc} Fuchs, B.A., Shabat, B.V. \emph{Functions of a Complex Variable and some of their Applications, Volume I.} Pergamon Press, 1964.

\bibitem{ham} Hammons Jr., A.R., Kumar, P.V., Calderbank, A.R., Sloane, N.J.A. \& Sol\'e, P. \emph{The $\zfour$-Linearity of Kerdock, Preparata, Goethals and Related Codes.} IEEE Trans Information Theory, 40, 1994, pp 301-319.

\bibitem{hed} Hedayat, A.S., Sloane, N.J.A. \& Stufken, J. \emph{Orthogonal Arrays: Theory and Applications.} Springer-Verlag, New York, 1999.

\bibitem{hun}
Hungerford, Thomas W., \textit{Algebra}. Springer-Verlag, New York, 1974.

\bibitem{mat} Matsumura, M., \textit{Commutative Ring Theory}, Cambridge
University Press, 1989.

\end{thebibliography}
\end{document}